\documentclass[11pt]{amsart}

\usepackage{amssymb}
\usepackage{latexsym}
\usepackage{amsmath}

\begin{document}

\newtheorem{theorem}{Theorem}[section]
\newtheorem{definition}{Definition}[section]
\newtheorem{proposition}{Proposition}[section]
\newtheorem{lemma}{Lemma}[section]

\title{Graphs having no quantum symmetry}

\author{T. Banica}
\address{T.B.: Department of Mathematics, Universite Toulouse 3, 118 route de Narbonne, 31062 Toulouse, France}
\email{banica@picard.ups-tlse.fr}
\author{J. Bichon}
\address{J.B.: Department of Mathematics, Universite
de Pau, 1 avenue de l'universite, 64000 Pau, France}
\email{bichon@univ-pau.fr}
\author{G. Chenevier}
\address{G.C.: Department of Mathematics, Universite
Paris 13, 99 avenue J-B. Clement, 93430 Villetaneuse, France}
\email{chenevie@math.univ-paris13.fr}


\begin{abstract}
We consider circulant graphs having $p$ vertices, with $p$ prime.
To any such graph we associate a certain number $k$, that we call
type of the graph. We prove that for $p>>k$ the graph has no
quantum symmetry, in the sense that the quantum automorphism group
reduces to the classical automorphism group.
\end{abstract}

\maketitle

\section*{Introduction}

A remarkable fact, discovered by Wang in \cite{wa}, is that the
set $\{1,\ldots ,n\}$ has a quantum permutation group. For
$n=1,2,3$ this the usual symmetric group $S_n$. However, starting
from $n=4$ the ``quantum permutations'' do exist. They form a
compact quantum group $\mathcal Q_n$, satisfying the axioms of
Woronowicz in \cite{wo}.

The next step is to look at ``simplest'' subgroups of $\mathcal
Q_n$. There are many natural degrees of complexity for such a
subgroup, and the notion that emerged is that of quantum
automorphism group of a vertex-transitive graph. These graphs are
those having the property that the usual automorphism group acts
transitively on the set of vertices. We assume of course that the
number of vertices is $n$.

These quantum groups are studied in \cite{bi1}, \cite{bi2} and
\cite{ba2}, \cite{ba3}, then in \cite{bb1}, \cite{bb2}.

The motivation comes from certain combinatorial aspects of
subfactors, free probability, and statistical mechanical models.
See \cite{ba3}, \cite{bb1}, \cite{bc}.

A fascinating question here, whose origins go back to Wang's paper
\cite{wa}, is to decide whether a given graph has quantum symmetry
or not. There are basically two kinds of graphs where the answer
is understood, namely: \begin{enumerate} \item The $n$-element set
$X_n$. This graph has $n$ vertices, and no edges at all. \item The
$n$-cycle $C_n$. This graph has $n$ vertices, $n$ edges, and looks
like a cycle.
\end{enumerate}

The graphs having no quantum symmetry are as follows:
\begin{enumerate}
\item $X_n$, $n<4$. This is proved in \cite{wa}, by direct algebraic computation. An explanation is proposed in \cite{ba1}, where the number $n\in{\mathbb N}$ is interpreted as a Jones index. This is further refined in \cite{ba3}, where
$\mathcal Q_n$ is shown to appear as Tannakian realisation of the Temperley-Lieb planar algebra of index $n$, known to be degenerate in the index range $1\leq n<4$.
\item $C_n$, $n\neq 4$. This is proved in \cite{ba2}, by direct algebraic computation. An explanation regarding $C_4$ is proposed in \cite{bb1}: this graph is exceptional in the series because it is the one having non-trivial disconnected complement. Indeed, the quantum symmetry group is the same for a graph and for its complement, and duplication of graphs corresponds to free wreath products, known from \cite{bi2} to be highly non-commutative operations.
\end{enumerate}

Some other results on lack of quantum symmetry include verifications for a number of cycles with chords, for a special graph called discrete torus, and stability/not stability under various product operations. See \cite{ba3}, \cite{bb1}, \cite{bb2}.

Although most such results have ad-hoc proofs, there is an idea emerging from this work, namely that computations become simpler with  $n\to\infty$.

In this paper we find an asymptotic result of non-quantum
symmetry. We consider graphs which are circulant, and have prime
number of vertices: that is, if $p$ is the number of vertices,
then $\mathbb Z_p$ must act on the graph, and $p$ must be prime.
To any such graph we associate a number $k$, that we call type,
and which measures in a certain sense the complexity of the graph
(as an example, for $C_n$ we have $k=2$). Our result is that a
type $k$ graph having enough vertices has no quantum symmetry.

The proof uses a standard technique, gradually developed since
Wang's paper \cite{wa}, and pushed here one step forward, by
combination with a Galois theory argument. We should mention that
the combination is done only at the end: it is not clear how to
include in the coaction formalism the underlying arithmetics.

We don't know what happens when the number of vertices is not prime:
\begin{enumerate}
\item Most ingredients have extensions to the general case, and it
won't be surprising that some kind of asymptotic result holds here
as well. However, there are a number of obstructions to be
overcome. These seem to come from complexity of the usual
automorphism group. For a prime number of vertices this group is
quite easy to describe, as shown by Alspach in \cite{al}, but in
general the situation is quite complicated, as shown for instance
by Klin and P\"oschel in \cite{kp}, or by Dobson and Morris in
\cite{dm}. \item A vertex-transitive graph having a prime number
of vertices is necessary circulant. So, in order to extend our
result, it is not clear whether to remain or not in the realm of
circulant graphs. Moreover, it would be interesting to switch at
some point to higher combinatorial structures, describing
arbitrary subgroups of $\mathcal Q_n$. In other words, there is a
lot of work to be done, and this paper should be regarded as a
first one on the subject.
\end{enumerate}

We should probably say a word about the original motivating
problems. As explained in \cite{ba2}, \cite{ba3}, \cite{bc},
quantum permutation groups are closely related to the ``$2$-box'',
``spin model'' and ``meander'' problems, discussed in \cite{bj},
\cite{df}, \cite{js}. We think that the idea in this paper is new
in the area -- for instance, it is not of topological nature --
and it is our hope that further developments of it, along the
above lines, might be of help in connection with these problems.

Finally, let us mention that the idea of letting $n\to\infty$ is
very familiar in certain areas of representation theory, developed
by Weingarten (\cite{wg}), Biane (\cite{bn}), Collins (\cite{co})
and many others. For quantum groups such methods are worked out in
\cite{bc}, but their relation with the present results is very
unclear.

The paper is organized as follows. Sections 1--2 are a quick
introduction to the problem, in 3 we fix some notations, and in
4--5 we prove the main result.

\subsection*{Acknowledgements} We would like to express our gratitude to the NLS research center in Paris and to the Institute for theoretical physics at Les Houches, for their warm hospitality and support, at an early stage of this project.

\section{Magic unitary matrices}

In this section and the next two ones we present a few basic facts
regarding quantum permutation groups, along with some examples,
explanations and sketches of proofs. The material is listed
according to an ad-hoc ordering, with the combinatorial side of
the subject coming first. For further reading, we recommend
\cite{ba3}.

Let $A$ be a ${\mathbb C}^*$-algebra. That is, we have a complex algebra with a norm and an involution, such that Cauchy sequences converge, and $||aa^*||=||a||^2$.

The basic examples are $B(H)$, the algebra of bounded operators on a Hilbert space $H$, and ${\mathbb C}(X)$, the algebra of continuous functions on a compact space $X$.

In fact, any ${\mathbb C}^*$-algebra is a subalgebra of some $B(H)$, and any commutative ${\mathbb C}^*$-algebra is of the form ${\mathbb C}(X)$. These are results of Gelfand-Naimark-Segal and Gelfand, both related to the spectral theorem for self-adjoint operators.

\begin{definition}
Let $A$ be a ${\mathbb C}^*$-algebra.
\begin{enumerate}
\item A projection is an element $p\in A$ satisfying $p^2=p=p^*$.
\item Two projections $p,q\in A$ are called orthogonal when
$pq=0$. \item A partition of unity is a set of orthogonal
projections, which sum up to $1$.
\end{enumerate}
\end{definition}

A projection in $B(H)$ is an orthogonal projection $\pi(K)$, where $K\subset H$ is a closed subspace. Orthogonality of projections corresponds to orthogonality of subspaces, and partitions of unity correspond to decompositions of $H$.

A projection in ${\mathbb C}(X)$ is a characteristic function $\chi(Y)$, where $Y\subset X$ is an open and closed subset. Orthogonality of projections corresponds to disjointness of subsets, and partitions of unity correspond to partitions of $X$.

\begin{definition}
A magic unitary is a square matrix $u\in M_n(A)$, all whose rows and columns are partitions of unity in $A$.
\end{definition}

Such a matrix is indeed unitary, in the sense that we have $uu^*=u^*u=1$.

Over $B(H)$ these are the matrices $\pi(K_{ij})$ with $K_{ij}$ magic decomposition of $H$, meaning that each row and column of $K$ is a decomposition of $H$.

Over ${\mathbb C}(X)$ these are the matrices $\chi(Y_{ij})$ with $Y_{ij}$ magic partition of $X$, meaning that each row and column of $Y$ is a partition of $X$.

We are interested in the following example.
Consider a finite graph $X$.
In this paper this means that we have a finite set of vertices,
and certain pairs of distinct vertices are connected by unoriented edges.
We do not allow multiple edges.

\begin{definition}
The magic unitary of a finite graph $X$ is given by $$u_{ij}=\chi\{g\in G\mid g(j)=i\}$$
where $i,j$ are vertices of $X$, and $G$ is the automorphism group of $X$.
\end{definition}

This is by definition a $V\times V$ matrix over the algebra $A={\mathbb C}(G)$, where $V$ is the vertex set. In case vertices are labeled $1,\ldots,n$, we can write $u\in M_n(A)$.

The fact that the characteristic functions $u_{ij}$ form indeed a magic unitary follows from the fact that the corresponding sets form a magic partition of $G$.

We denote by $d$ the adjacency matrix of $X$. This is a $V\times
V$ matrix, given by $d_{ij}=1$ if $i,j$ are connected by an edge,
and by $d_{ij}=0$ if not.

We have the following presentation result.

\begin{theorem}
The algebra $A={\mathbb C}(G)$ is isomorphic to the universal ${\mathbb C}^*$-algebra generated by $n^2$ elements $u_{ij}$, with the following relations:
\begin{enumerate}
\item The matrix $u=(u_{ij})$ is a magic unitary.
\item We have $du=ud$, where $d$ is the adjacency matrix of $X$.
\item The elements $u_{ij}$ commute with each other.
\end{enumerate}
\end{theorem}

\begin{proof}
Let $A'$ be the universal algebra in the statement. That is, $A'$
is the universal repelling object in the category of commutative
$\mathbb C^*$-algebras generated by entries of a $n\times n$ magic
unitary matrix $u$, subject to the condition $du=ud$. The
construction of such an object is standard, and we have uniqueness
up to isomorphism.

The magic unitary of $X$ commutes with $d$, so we have a morphism
$p:A'\to A$. By applying Gelfand's theorem, $p$ comes from an
inclusion $i:G\subset G'$, where $G'$ is the spectrum of $A'$.

By using the universal property of $A'$, we see that the formulae
\begin{eqnarray*}
\Delta(u_{ij})&=&\sum u_{ik}\otimes u_{kj}\cr
\varepsilon (u_{ij})&=&\delta_{ij}\cr
S(u_{ij})&=&u_{ji}
\end{eqnarray*}
define morphisms of algebras. These must come from maps $G'\times
G',\{.\},G'\to G'$, making $G'$ into a group, acting on $X$, and
we get $G=G'$. See section 2 in \cite{ba3} for missing details.
\end{proof}

\section{Quantum permutation groups}

Let $X$ be a graph as in previous section. Its quantum automorphism group is constructed by removing commutativity from Theorem 1.1 and its proof.

\begin{definition}
The Hopf algebra associated to $X$ is the universal ${\mathbb C}^*$-algebra ${\mathcal A}$ generated by entries $u_{ij}$ of a $n\times n$ magic unitary commuting with $d$, with
\begin{eqnarray*}
\Delta(u_{ij})&=&\sum u_{ik}\otimes u_{kj}\cr
\varepsilon (u_{ij})&=&\delta_{ij}\cr
S(u_{ij})&=&u_{ji}
\end{eqnarray*}
as comultiplication, counit and antipode maps.
\end{definition}

The precise structure of ${\mathcal A}$ is that of a co-involutive
unital Hopf ${\mathbb C}^*$-algebra of finite type. That is,
${\mathcal A}$ satisfies the axioms of Woronowicz in \cite{wo},
along with the extra axiom $S^2=id$. See \cite{ba3}, \cite{ks} for
more details on this subject.

For the purposes of this paper, let us just mention that we have the formula
$${\mathcal A}={\mathbb C}({\mathcal G})$$
where ${\mathcal G}$ is a compact quantum group. This quantum group doesn't exist as a concrete object, but several tools from Woronowicz's paper \cite{wo}, such as an analogue of the Peter-Weyl theory, are available for it, in the form of functional analytic statements  regarding its algebra of continuous functions ${\mathcal A}$.

Comparison of Theorem 1.1 and Definition 2.1 shows that we have a morphism ${\mathcal A}\to {\mathbb C}(G)$. This can be thought of as coming from an inclusion $G\subset{\mathcal G}$.

\begin{definition}
We say that $X$ has no quantum symmetry if ${\mathcal A}={\mathbb C}(G)$.
\end{definition}

It is not clear at this point whether there exist graphs $X$ which do have quantum symmetry. Before getting into the subject, let us state the following useful result.

\begin{theorem}
The following are equivalent.
\begin{enumerate}
\item $X$ has no quantum symmetry.
\item ${\mathcal A}$ is commutative.
\item For $u$ magic unitary, $du=ud$ implies that $u_{ij}$ commute with each other.
\end{enumerate}
\end{theorem}

\begin{proof}
All equivalences are clear from definitions, and from the Gelfand theorem argument in proof of Theorem 1.1.
\end{proof}

The very first graphs to be investigated are the $n$-element sets $X_n$. Here the incidency matrix is $d=0$, so the above condition (3) is that for any $n\times n$ magic unitary matrix $u$, the entries $u_{ij}$ have to commute with each other.
\begin{enumerate}
\item The graph $X_2$. This has no quantum symmetry, because a $2\times 2$ magic unitary has to be of the form
$$u_p=\begin{pmatrix}
p&1-p\cr 1-p&p\end{pmatrix}$$
with $p$ projection, and entries of this matrix commute with each other.
\item The graph $X_3$. This has no quantum symmetry either, as shown in \cite{wa}.
\item The graph $X_4$. This has quantum symmetry, because the matrix
$$u_{pq}=\begin{pmatrix}
p&1-p&0&0\cr
1-p&p&0&0\cr
0&0&q&1-q\cr
0&0&1-q&q
\end{pmatrix}$$
is a magic unitary, whose entries don't commute if $pq\neq qp$.
\item The graph $X_{n}$, $n\geq 5$. This has quantum symmetry too,
as one can see by adding to $u_{pq}$ a diagonal tail formed of
$1$'s.
\end{enumerate}

The other series of graphs where complete results are available
are the $n$-cycles $C_n$. The situation here, already described in
the introduction, is as follows.
\begin{enumerate}
\item The graph $C_2$. This has no quantum symmetry, because $X_2$ doesn't.
\item The graph $C_3$. This has no quantum symmetry, because $X_3$ doesn't.
\item The graph $C_4$. This has quantum symmetry, because its adjacency matrix
$$d=\begin{pmatrix}
0&0&1&1\cr
0&0&1&1\cr
1&1&0&0\cr
1&1&0&0
\end{pmatrix}$$
written here according to the assignment of numbers $1324$ to
vertices in a cyclic way, commutes with $u_{pq}$. \item The graph
$C_n$, $n\geq 5$. This has no quantum symmetry, as shown in
\cite{ba2}.
\end{enumerate}

Summarizing, the subtle results in these series are those
regarding lack of quantum symmetry of cycles $C_n$, with $n=3$ and
$n\geq 5$. In what follows we present a general result, which
applies in particular to $C_p$ with $p$ big prime (in fact $p\geq
7$). This result will have the following consequences to what has
been said so far:
\begin{enumerate}
\item As explained in the introduction, we hope to extend at some
point our techniques, as to apply to $C_n$ with big $n$. \item As
for $C_n$ with small $n$, we won't think about it for some time.
This is an exceptional graph, at least until the asymptotic area
is well understood.
\end{enumerate}

Our last remark is about use of $\mathbb C^*$-algebras. The lack
of quantum symmetry can be characterized in fact in a purely
algebraic manner. Indeed, consider $\mathcal A_0$, the universal
$*$-algebra generated by entries $u_{ij}$ of a $n\times n$ magic
unitary matrix commuting with $d$. By using general theory from
\cite{ks}, namely Theorem 27 and Proposition 32 in Chapter 11, we
get a $*$-algebra embedding with dense image $\mathcal A_0 \to \mathcal A$. This
shows that $\mathcal A$ is commutative if and only if $\mathcal
A_0$ is.

\section{Circulant graphs}

A graph $X$ having $n$ vertices is called circulant if its automorphism group contains a cycle of length $n$, and hence a copy of the cyclic group ${\mathbb Z}_n$.

This is the same as saying that vertices of $X$ are $n$-th roots of unity, edges are represented by certain segments, and the whole picture has the property of being invariant under the $2\pi/n$ rotation centered at $0$. Here the rotation is either the clockwise or the counterclockwise one: the two conditions are equivalent.

For the purposes of this paper, best is to assume that vertices of
$X$ are elements of ${\mathbb Z}_n$, and $i\sim j$ (connection by
an edge) implies $i+k\sim j+k$ for any $k$.

We denote by ${\mathbb Z}_n^*$ the group of invertible elements of
the ring ${\mathbb Z}_n$.

Our study of circulant graphs is based on diagonalisation of corresponding adjacency matrices. This is in turn related to certain arithmetic invariants of the graph -- an abelian group $E$ and a number $k$ -- constructed in the following way.

\begin{definition}
Let $X$ be a circulant graph on $n$ vertices.
\begin{enumerate}
\item The set $S\subset {\mathbb Z}_n$ is given by $i\sim j \iff j-i \in S$.
\item The group $E\subset{\mathbb Z}_n^*$ consists of elements $a$ such that $aS=S$.
\item The order of $E$ is denoted $k$, and is called type of $X$.
\end{enumerate}
\end{definition}

The interest in $k$ is that this is the good parameter measuring complexity of the spectral theory of $X$. Calling it ``type'' might seem a bit unnatural at this point; but the terminology will be justified by the main result in this paper.

Here are a few basic examples and properties, $\varphi$ being the
Euler function:
\begin{enumerate}
\item The type can be $2,4,6,8,\ldots$ This is because $\{\pm 1\}\subset E$.
\item $C_n$ is of type $2$. Indeed, we have $S=\{\pm 1\}$, $E=\{\pm 1\}$.
\item $X_n$ is of type $\varphi(n)$. Indeed, here $S=\emptyset$, $E={\mathbb Z}_n^*$.
\end{enumerate}

It is possible to make an extensive study of this notion, but we
won't get into the subject. Let us just mention that the graphs
$2C_5,C_{10}$ studied in \cite{bb2} have the same $E$ group, but
the first one has quantum symmetry, while the second one hasn't.
Here $2C_5$ is the disjoint union of two copies of $C_5$, and the
fact that this graph has quantum symmetry comes from free wreath
product philosophy.

Consider the Hopf algebra ${\mathcal A}$ associated to $X$, as in previous section.

\begin{definition}
The linear map $\alpha :{\mathbb C}^n\to{\mathbb C}^n\otimes{\mathcal A}$ given by the formula
$$\alpha (e_i)=\sum e_j\otimes u_{ji}$$
where $e_1,\ldots ,e_n$ is the canonical basis of ${\mathbb C}^n$, is called coaction of ${\mathcal A}$.
\end{definition}

It follows from the magic unitarity condition that $\alpha$ is a morphism of algebras, which satisfies indeed the axioms of coactions. See \cite{ba3} for details.

For the purposes of this paper, let us just mention that $\alpha$
appears as functional analytic transpose of the action of
${\mathcal G}$ on the set $X_n=\{1,\ldots ,n\}$. In other words,
we have $\alpha (\varphi)=\varphi\circ a$, where $a:X_n\times
{\mathcal G}\to X_n$ is the action map, $a(i,g)=g(i)$.

These general considerations are valid in fact for any graph. In  what follows we use the following simple fact, valid as well in the general case.

\begin{theorem}
If $F$ is an eigenspace of $d$ then $\alpha(F)\subset F\otimes {\mathcal A}$.
\end{theorem}

\begin{proof}
Since $u$ commutes with $d$, it commutes with the ${\mathbb
C}^*$-algebra generated by $d$, and in particular with the
projection $\pi(F)$. The relation $u\pi(F)=\pi(F)u$ can be
translated in terms of $\alpha$, and we get $\alpha(F)\subset
F\otimes {\mathcal A}$. See section 2 in \cite{ba3}.
\end{proof}

\section{Spectral decomposition}

In what follows $X$ is a circulant graph having $p$ vertices, with $p$ prime.

We denote by $d,{\mathcal A},\alpha$ the associated adjacency matrix, Hopf algebra and coaction, and by $S,E,k$ the set, group and number in Definition 3.1.

We denote by $\xi$ the column vector $(1,w,w^2,\ldots ,w^{p-1})$, where $w=e^{2\pi i/p}$.

\begin{lemma}
The eigenspaces of $d$ are given by $V_0={\mathbb C}1$ and
$$V_x=\bigoplus_{a\in E}{\mathbb C}\,\xi^{xa}$$
with $x\in \mathbb Z_p^*$. Moreover, we have $V_x=V_y$ if and only if $xE=yE$.
\end{lemma}

\begin{proof}
The matrix $d$ being circulant, we have the formula
$$d(\xi^x)=f(x)\xi^x$$
where $f:{\mathbb Z}_p\to{\mathbb C}$ is the following function:
$$f(x)=\sum_{t\in S}w^{xt}$$

Let $K={\mathbb Q}(w)$ and let $H$ be the Galois group of the Galois extension $\mathbb Q \subset K$.
It is well-known that we have a group isomorphism \begin{align*}
\mathbb Z_p^* & \longrightarrow H \\
x & \longmapsto s_x
\end{align*}
with the automorphism $s_x$ given by the following formula:
$$s_x({w}) = {w}^x$$

Also, we know from a theorem of Dedekind that the family $\{s_x\mid x\in{\mathbb Z}_p^*\}$ is free in ${\rm End}_{\mathbb Q}(K)$. Now for $x,y\in \mathbb Z_p^*$ consider the following operator:
$$L = \sum_{t \in S} s_{xt} - \sum_{t \in S} s_{yt} \in
{\rm End}_{\mathbb Q}(K)$$

We have $L({w}) = f(x)-f(y)$, and since $L$ commutes with the action of the
abelian group $H$, we have
$$L=0 \iff L({w}) =0 \iff f(x)=f(y)$$
and by linear independence of the family $\{s_x\mid x\in \mathbb Z_p^*\}$ we get:
$$f(x) = f(y) \iff xS=yS \iff xE=yE$$

It follows that $d$ has precisely $1+(p-1)/k$ distinct eigenvalues, the corresponding
eigenspaces being those in the statement.
\end{proof}

Consider now a commutative ring $(R,+,\cdot)$. We denote by $R^*$
the group of invertibles, and we assume $2\in R^*$. A subgroup
${G}\subset R^*$ is called even if $-1\in{G}$.

\begin{definition}
An even subgroup ${G}\subset R^*$ is called $2$-maximal if
$$a-b=2(c-d)$$
with $a,b,c,d\in{G}$ implies $a=\pm b$.
\end{definition}

We call $a=b,c=d$ trivial solutions, and $a=-b=c-d$ hexagonal
solutions. The terminology comes from the following key example:

Consider the group $G\subset{\mathbb C}$ formed by $k$-th roots of
unity, with $k$ even. We regard $G$ as set of vertices of the
regular $k$-gon. An equation of the form $a-b=2(c-d)$ with
$a,b,c,d\in G$ says that the diagonals $a-b$ and $c-d$ are
parallel, and that the first one is twice as much as the second
one. But this can happen only when $a,c,d,b$ are consecutive
vertices of a regular hexagon, and here we have $a+b=0$.

This example is discussed in detail in next section.

\begin{proposition}
Assume that $R$ has the property $3\neq 0$, and consider a
$2$-maximal subgroup ${G}\subset R^*$.
\begin{enumerate}
\item $2,3\not\in{G}$. \item $a+b=2c$ with $a,b,c\in{G}$ implies
$a=b=c$. \item $a+2b=3c$ with $a,b,c\in{G}$ implies $a=b=c$.
\end{enumerate}
\end{proposition}

\begin{proof}
(1) This follows from the following formulae, which cannot hold in
${G}$:
$$4-2=2(2-1)$$
$$3-(-1)=2(3-1)$$

Indeed, the first one would imply $4=\pm 2$, and the second one
would imply $3=\pm 1$. But from $2\in R^*$ and $3\neq 0$ we get
$2,4,6\neq 0$, contradiction.

(2) We have $a-b=2(c-b)$. For a trivial solution we have
 $a=b=c$, and for a hexagonal
 solution we have $a+b=0$,
hence $c=0$, hence $0\in{G}$, contradiction.

(3) We have $a-c=2(c-b)$. For a trivial solution we have
 $a=b=c$, and for a hexagonal
 solution we have $a+c=0$,
hence $b=-2a$, hence $2\in{G}$, contradiction.
\end{proof}

We use these facts several times in the proof below, by refering
to them as ``$2$-maximality'' properties, without special mention
to Proposition 4.1.

\begin{theorem}
If $E\subset {\mathbb Z}_p$ is $2$-maximal ($p\geq 5$)
then $X$ has no quantum
symmetry.
\end{theorem}

\begin{proof}
We use Lemma 4.1, which ensures that $V_1,V_2,V_3$ are eigenspaces of $d$. By $2$-maximality of $E$, these
three eigenspaces are different.

From eigenspace preservation in Theorem 3.1 we get formulae of the
following type, with $r_a,r_a',r_a''\in{\mathcal A}$:
\begin{eqnarray*}
\alpha (\xi)&=&\sum_{a\in E}\xi^a\otimes r_a\cr
\alpha (\xi^2)&=&\sum_{a\in E}\xi^{2a}\otimes r_a'\cr
\alpha (\xi^3)&=&\sum_{a\in E}\xi^{3a}\otimes r_a''
\end{eqnarray*}

We take the square of the first relation, we compare with the
formula of $\alpha(\xi^2)$, and we use $2$-maximality:
\begin{eqnarray*}
\alpha(\xi^2) &=&\left( \sum_{a\in E}\xi^a\otimes r_a\right)^2\cr
&=&\sum_{x}\xi^x\otimes \left(\sum_{a,b\in
E}\delta_{a+b,x}\,r_ar_b\right)\cr &=&\sum_{c\in E}\xi^{2c}\otimes
\left(\sum_{a,b\in E}\delta_{a+b,2c}\,r_ar_b\right)\cr
&=&\sum_{c\in E}\xi^{2c}\otimes r_c^2
\end{eqnarray*}

We multiply this relation by the formula of $\alpha(\xi)$, we
compare with the formula of $\alpha(\xi^3)$, and we use
$2$-maximality:
\begin{eqnarray*}
\alpha(\xi^3) &=&\left(\sum_{a\in E}\xi^a\otimes r_a\right)
\left(\sum_{c\in E}\xi^{2c}\otimes r_c^2\right)\cr
&=&\sum_{x}\xi^x\otimes \left(\sum_{a,c\in
E}\delta_{a+2c,x}\,r_ar_c^2\right)\cr &=&\sum_{b\in
E}\xi^{3b}\otimes \left(\sum_{a,c\in
E}\delta_{a+2c,3b}\,r_ar_c^2\right)\cr &=&\sum_{b\in
E}\xi^{3b}\otimes r_b^3
\end{eqnarray*}

Summarizing, the three formulae in the beginning are in fact:
\begin{eqnarray*}
\alpha (\xi)&=&\sum_{a\in E}\xi^a\otimes r_a\cr \alpha
(\xi^2)&=&\sum_{a\in E}\xi^{2a}\otimes r_a^2\cr \alpha
(\xi^3)&=&\sum_{a\in E}\xi^{3a}\otimes r_a^3
\end{eqnarray*}

We claim now that for $a\neq b$, we have the following ``key
formula'':
$$r_ar_b^3=0$$

Indeed, consider the following equality:
$$\left(\sum_{a\in E}\xi^a\otimes r_a\right)
\left(\sum_{b\in E}\xi^{2b}\otimes r_b^2\right)=\sum_{c\in
E}\xi^{3c}\otimes r_c^3$$

By eliminating all $a=b$ terms, which produce the sum on the
right, we get:
$$\sum\left\{\xi^{a+2b}\otimes r_ar_b^2\mid a,b\in E,\,a
\neq b\right\}=0$$

By taking the coefficient of $\xi^x$, with $x$ arbitrary, we get:
$$\sum\left\{r_ar_b^2\mid a,b\in E,\,a\neq b,\,a+2b=x\right\}=0$$

We fix now $a,b\in E$ satisfying $a\neq b$. We know from
$2$-maximality that the equation $a+2b=a'+2b'$ with $a',b'\in E$
has at most one non-trivial solution, namely the hexagonal one,
given by $a'=-a$ and $b'=a+b$. Now with $x=a+2b$, we get that the
above equality is in fact one of the following two equalities:
$$r_ar_b^2=0$$
$$r_ar_b^2+r_{-a}r_{a+b}^2=0$$

In the first situation, we have $r_ar_b^3=0$ as claimed.

In the second situation, we proceed as follows. We know that
$a_1=b$ and $b_1=a+b$ are distinct elements of $E$. Consider now
the equation $a_1+2b_1=a_1'+2b_1'$ with $a_1',b_1'\in E$. The
hexagonal solution of this equation, given by $a_1'=-a_1$ and
$b_1'=a_1+b_1$, cannot appear: indeed, $b_1'=a_1+b_1$ can be
written as $b_1'=a+2b$, and by $2$-maximality we get $b_1'=-a=b$,
which contradicts $a+b\in E$.

Thus the equation $a_1+2b_1=a_1'+2b_1'$ with $a_1',b_1'\in E$ has
only trivial solutions, and with $x=a_1+2b_1$ in the above
considerations we get:
$$r_{a_1}r_{b_1}^2=0$$

Now remember that this follows by identifying coefficients in
$\alpha(\xi)\alpha(\xi^2)=\alpha(\xi^3)$. The same method applies
to the formula $\alpha(\xi^2)\alpha(\xi)=\alpha(\xi^3)$, and we
get:
$$r_{b_1}^2r_{a_1}=0$$

We have now all ingredients for finishing the proof of the key
formula:
\begin{eqnarray*}
r_ar_b^3 &=&r_ar_b^2r_b\cr &=&-r_{-a}r_{a+b}^2r_b\cr
&=&-r_{-a}r_{b_1}^2r_{a_1}\cr&=&0
\end{eqnarray*}

We come back to the following formula, proved for $s=1,2,3$:
$$\alpha(\xi^s)=\sum_{a\in E}\xi^{sa}\otimes r_a^s$$

By using the key formula, we get by induction on $s\geq 3$ that
this holds in general:
\begin{eqnarray*}
\alpha\left(\xi^{1+s}\right) &=&\left(\sum_{a\in E}\xi^a\otimes
r_a\right) \left(\sum_{b\in E}\xi^{sb}\otimes r_b^s\right)\cr
&=&\sum_{a\in E}\xi^{(1+s)a}\otimes r_a^{1+s}+\sum_{a,b\in
E,\,a\neq b}\xi^{a+sb}\otimes r_ar_b^s\cr &=&\sum_{a\in
E}\xi^{(1+s)a}\otimes r_a^{1+s}
\end{eqnarray*}

In particular with $s=p-1$ we get:
$$\alpha(\xi^{-1})=\sum_{a\in E}\xi^{-a}\otimes r_a^{p-1}$$

On the other hand, from $\xi^*=\xi^{-1}$ we get
$$\alpha(\xi^{-1})=\sum_{a\in E}\xi^{-a}\otimes r_a^*$$
which gives $r_a^* = r_a^{p-1}$ for any $a$. Now by using the key
formula we get
$$(r_ar_b)(r_ar_b)^*=r_ar_br_b^*r_a^*=r_ar_b^pr_a^*=(r_ar_b^3)(r_b^{p-3}r_a^*)=0$$
which gives $r_ar_b=0$. Thus we have $r_ar_b=r_br_a=0$.

On the other hand, ${\mathcal A}$ is generated by coefficients of
$\alpha$, which are in turn powers of elements $r_a$. It follows
that ${\mathcal A}$ is commutative, and we are done.
\end{proof}

\section{The main result}

Let $k$ be an even number, and consider the group of $k$-th roots
of unity $G=\{1,\zeta,\ldots ,\zeta^{k-1}\}$, where $\zeta=e^{2\pi
i/k}$. We use the Euler function $\varphi$.

\begin{lemma}
$G$ is $2$-maximal in ${\mathbb C}$.
\end{lemma}

\begin{proof}
Assume that we have $a-b=2(c-d)$ with $a,b,c,d\in G$. With $z=b/a$
and $u=(c-d)/a$, we have $1-z=2u$. Let $n$ be the order of the
root of unity $z$. By \cite{ws}, chap. 2, the $\mathbb Q(z)$-norm
$N(1-z)$ of $1-z$ is $\pm 1$ if $n$ is not the power of a prime
$l$, and $\pm l$ otherwise. Applying the $\mathbb Q(z)$-norm to
$1-z=2u$, and using that $u$ is an algebraic integer, we get
$$2^{\varphi(n)}\, \vert \, N(1-z)$$ hence $n\leq 2$, $z=\pm 1$,
and we are done.
\end{proof}

Let $p$ be a prime number.

\begin{lemma}
For $p>6^{\varphi(k)}$, any subgroup $E\subset{\mathbb Z}_p^*$ of
order $k$ is $2$-maximal.
\end{lemma}

\begin{proof}
Consider the following set of complex numbers:
$$\Sigma=\{a+2b\mid a,b\in G\}$$

Let $A={\mathbb Z}[\zeta]$, recall that $A$ is the ring of
algebraic integers of ${\mathbb Q}(\zeta)$, and in particular a
Dedekind ring. If $p$ is any prime number such that $k$ divides
$p-1$, it is well-known that the ideal $pA$ is a product
$P_{1}\ldots P_{\varphi(k)}$ of prime ideals of $A$ such that
$A/P_{i} \simeq {\mathbb Z}_{p}$ for each $i$. Choosing an $i$ we
get a surjective ring morphism:
$$\Phi: A \rightarrow {\mathbb Z}_{p}$$

Since $p$ does not divide $k$, the polynomial
$$X^{k}-1=\prod_{i=0}^{k-1} \left(X-\Phi(\zeta)^{i}\right)$$
has no multiple root in ${\mathbb Z}_{p}$, hence $\Phi(G) \subset
{\mathbb Z}_{p}^{*}$ is a cyclic subgroup of order $k$. As
${\mathbb Z}_{p}^{*}$ is known to be a cyclic group, $\Phi(G)$ is
actually the unique subgroup of order $k$ of ${\mathbb
Z}_{p}^{*}$, hence it coincides with the subgroup $E$ in the
statement.

We claim that for $p$ as in the statement, the induced map $\Phi:
\Sigma \rightarrow {\mathbb Z}_{p}$ is injective. Together with
Lemma 5.1, this would prove the assertion.

So, assume $\Phi(x)=\Phi(y)$. The Dedekind property gives an ideal
$Q\subset A$ such that:
$$ (x-y)=P_{i}Q$$

For $I$ a nonzero ideal of $A$, let us denote by $N(I):=\big\vert
A/I\big\vert$ the norm of $I$, and set also $N(0)=0$. Recall that
by the Dedekind property, $N$ is multiplicative with respect to
the product of ideals in $A$ and that for any $z\in A$, the norm
$N(z)$ of the principal ideal $zA$ coincides with the absolute
value of the following integer:
$$\prod_{s \in {\rm Gal}({\mathbb Q}(\zeta)/{\mathbb Q})}s(z)$$

Applying norms to $(x-y)=P_{i}Q$ shows that $N(P_{i})=p$ divides
the integer $N(x-y)$. Now with $p$ as in the statement we have
$N(x-y)\leq p_0$ for any $x,y\in \Sigma$, so the induced map
$\Phi:\Sigma\to{\mathbb Z}_p$ is injective, and we are done.
\end{proof}

\begin{theorem}
A type $k$ circulant graph having $p>>k$ vertices, with $p$ prime,
has no quantum symmetry.
\end{theorem}

\begin{proof}
This follows from Theorem 4.1 and Lemma 5.2, with
$p>6^{{\varphi(k)}}$.
\end{proof}


\begin{thebibliography}{99}
\bibitem{al}B. Alspach, Point-symmetric graphs and digraphs of prime order and transitive permutation groups of prime degree, {\em J. Combinatorial Theory Ser. B} {\bf 15} (1973), 12--17.
\bibitem{ba1}T. Banica, Symmetries of a generic coaction, {\em Math. Ann.} {\bf 314} (1999), 763--780.
\bibitem{ba2} T. Banica, Quantum automorphism groups of small metric spaces, {\em Pacific J. Math.} {\bf 219} (2005), 27--51.
\bibitem{ba3} T. Banica, Quantum automorphism groups of homogeneous graphs, {\em J. Funct. Anal.} {\bf 224} (2005), 243-280.
\bibitem{bb1}T. Banica and J. Bichon, Free product formulae for quantum permutation groups, {\em J. Math. Inst. Jussieu}, to appear.
\bibitem{bb2}T. Banica and J. Bichon, Quantum
automorphism groups of vertex-transitive graphs of order $\leq$
11, {\tt math.QA/0601758}.
\bibitem{bc}T. Banica and B. Collins, Integration over compact quantum groups,  {\em Publ. Res. Inst. Math. Sci.}, to appear.
\bibitem{bn}P. Biane, Representations of symmetric groups and free probability, {\em Adv. Math.} {\bf 138} (1998), 126--181.
\bibitem{bi1} J. Bichon, Quantum automorphism groups of finite graphs, {\em Proc. Amer. Math. Soc.} {\bf 131} (2003), 665--673.
\bibitem{bi2} J. Bichon, Free wreath product by the quantum permutation group, {\em Alg. Rep. Theory} {\bf 7} (2004), 343--362.
\bibitem{bj} D. Bisch and V.F.R. Jones, Singly generated planar algebras of small dimension, {\em Duke Math. J.} {\bf 101} (2000), 41--75.
\bibitem{co}B. Collins, Moments and cumulants of polynomial random variables on unitary groups, the {I}tzykson-{Z}uber integral, and free probability, {\em Int. Math. Res. Not.} {\bf 17} (2003), 953--982.
\bibitem{df} P. Di Francesco, Meander determinants,
{\em Comm. Math. Phys.} {\bf 191} (1998), 543--583.
\bibitem{dm}E. Dobson and J. Morris, On automorphism groups of circulant digraphs of square-free order, {\em Discrete Math.} {\bf 299} (2005), 79--98.
\bibitem{js} V.F.R. Jones and V.S. Sunder, Introduction to subfactors, LMS Lecture Notes 234, Cambridge University Press (1997)
\bibitem{ks}A. Klimyk and K. Schm\"udgen, Quantum groups and their representations, Texts and Monographs in Physics, Springer-Verlag, Berlin (1997).
\bibitem{kp}M.H. Klin and R. P\"oschel, The K\"onig problem, the isomorphism problem for cyclic graphs and the method of Schur rings, {\em Colloq. Math. Soc. Janos Bolyai} {\bf 25} (1981), 405--434.
\bibitem{wa} S. Wang, Quantum symmetry groups of finite spaces, {\em Comm. Math. Phys.} {\bf 195} (1998), 195--211.
\bibitem{ws}L.C. Washington, Introduction to cyclotomic fields, GTM 83, Springer (1982)
\bibitem{wg}D. Weingarten, Asymptotic behavior of group integrals in the limit of infinite rank, {\em J. Math. Phys.} {\bf 19} (1978), 999--1001.
\bibitem{wo}S.L. Woronowicz, Compact matrix pseudogroups, {\em Comm. Math. Phys.} {\bf 111} (1987), 613--665.
\end{thebibliography}
\end{document}